\documentclass{article}
\usepackage[utf8]{inputenc}
\usepackage{amsmath}
\usepackage{graphicx}
\usepackage{geometry}
\usepackage{hyperref}
\title{From Controllability to Information: A Unified Analysis via Gramian, Minimum Energy, Fisher Information Matrix and Entropy}
\author{Gabriel R. de Andrade Silva}
\date{July-01-2025}

\begin{document}

\maketitle

\begin{abstract}
This article explores the connections between controllability, control energy, information, and entropy in the context of linear dynamical systems, using the damped harmonic oscillator as a case study. We analytically derive the Controllability Gramian Matrix ($W_c$) and its determinant ($\det(W_c)$) for different damping regimes. We analyze the physical interpretation of $\det(W_c)$ in terms of minimum control energy, highlighting the influence of system parameters (damping $\zeta$ and natural frequency $\omega_n$). We investigate the conceptual and formal relationships between the Gramian, the Fisher Information Matrix ($I$), and Shannon ($H$) and thermodynamic ($S$) entropies, suggesting a duality between energetic control and estimation precision, as well as a link between ease of control and informational/thermodynamic entropy. The unified analysis offers insights into the fundamental principles governing the control of dynamical systems, integrating energy and information perspectives.
\end{abstract}

\section{Introduction}

Understanding and controlling complex dynamical systems represent central challenges in various areas of science and engineering, from physics and biology to economics and social sciences. The ability to predict and influence the behavior of these systems is fundamental for the development of advanced technologies, process optimization, and solving global problems. At the core of control theory, the notion of controllability, introduced by Kalman in the 1960s, establishes the conditions under which it is possible to drive a system from an initial state to any desired final state in finite time, using a limited set of control inputs. Although fundamental, classical controllability is a qualitative property, providing no information about the cost or effort required to achieve such a transition.

In practice, energy limitations are ubiquitous. Physical actuators have constraints on power, amplitude, and frequency, making purely qualitative analysis of controllability insufficient for many real applications. Thus arises the need for quantitative metrics that capture the control effort, allowing assessment of the ease or difficulty of maneuvering a system's state. The Controllability Gramian Matrix ($W_c$), defined as the integral of the exponential of the system's dynamics matrix multiplied by the input matrix and its transpose, has emerged as a powerful tool in this context. The spectral properties of the Gramian, such as its minimum and maximum eigenvalues, trace, and determinant, provide quantitative measures of the energy effort associated with control.

In particular, the determinant of the Controllability Gramian ($\det(W_c)$) gained prominence after the seminal work of Joshi and Mitter (1977), which established its fundamental physical interpretation: the inverse of the determinant is related to the minimum energy required to control the system. The higher the determinant, the lower the required minimum energy, indicating a system that is "easier" to control from an energetic perspective. This connection between an abstract mathematical property (the determinant of an integral matrix) and a measurable physical quantity (energy) opened new avenues for the analysis and design of efficient control systems.

In parallel, concepts from information theory and statistical physics, such as the Fisher Information Matrix ($I$) and Entropy (both in Shannon's formulation and thermodynamics), offer complementary perspectives for quantifying the uncertainty, information, and complexity inherent in dynamical systems. The Fisher Matrix, classically used in statistics to measure the amount of information that observable data carry about unknown parameters of a model, finds intriguing parallels in the control context, where it can be related to the curvature of the "reachable state space" and the distribution of control capability. Entropy, in turn, quantifies disorder or lack of information about a system's state. Recent investigations have explored the profound connections between these different metrics, suggesting that control energy, Fisher information, and entropy can be viewed as distinct facets of a more unified description of the dynamics and controllability of complex systems (Baggio et al., 2022).

This article seeks to explore these connections in depth, using the damped harmonic oscillator as a canonical model system. Through detailed analysis of this system, we will derive analytical expressions for the Controllability Gramian Matrix and its determinant in different damping regimes (underdamped, critically damped, overdamped, and undamped). We will investigate the physical interpretation of the determinant in terms of minimum control energy, explicitly showing the influence of damping on the effort required to control the oscillator. Next, we will explore the relationships between the Gramian, the Fisher Information Matrix, and Shannon and thermodynamic entropies, seeking to elucidate how these different metrics complement each other in characterizing the controllability and complexity of the system. The ultimate goal is to provide a unified analysis that integrates the concepts of controllability, energy, and information, offering insights into the fundamental principles governing the control of dynamical systems.

The work is organized as follows: Section 2 presents the mathematical modeling of the damped harmonic oscillator, including its state-space representation. Section 3 details the controllability analysis via the Gramian Matrix, with analytical calculation of the matrix and the matrix exponential for each regime. Section 4 focuses on the Gramian determinant, its analytical calculation, and its physical interpretation in terms of minimum control energy. Section 5 introduces the Fisher Information Matrix and explores its connection with the Gramian. Section 6 investigates the relationships with Shannon entropy and thermodynamic entropy. Section 7 presents a discussion of the results and their implications. Finally, Section 8 concludes the article with a recapitulation of the main contributions.

\section{Mathematical Modeling of the Harmonic Oscillator}

The damped harmonic oscillator is a fundamental and ubiquitous physical system, serving as a model for a vast range of phenomena in physics, engineering, and other sciences. Its dynamics are described by a second-order linear differential equation that captures the interaction between restoring forces (proportional to displacement), dissipative forces (proportional to velocity), and potentially external forces.

The general differential equation for the displacement $x(t)$ of a damped and forced harmonic oscillator is given by:

\begin{equation}
    m \ddot{x}(t) + c \dot{x}(t) + k x(t) = F(t)
\end{equation}

Where:
\begin{itemize}
    \item $m$ is the mass of the oscillator.
    \item $c$ is the viscous damping coefficient.
    \item $k$ is the spring constant (stiffness).
    \item $x(t)$ is the displacement of the oscillator relative to its equilibrium position at time $t$.
    \item $\dot{x}(t)$ is the velocity of the oscillator at time $t$.
    \item $\ddot{x}(t)$ is the acceleration of the oscillator at time $t$.
    \item $F(t)$ is the external force applied to the system at time $t$.
\end{itemize}

For systems and control analysis purposes, it is convenient to normalize this equation. Dividing by $m$ and introducing the parameters of undamped natural frequency ($\omega_n$) and damping factor ($\zeta$), we have:

\begin{equation}
    \ddot{x}(t) + \frac{c}{m} \dot{x}(t) + \frac{k}{m} x(t) = \frac{F(t)}{m}
\end{equation}

Defining:
\begin{itemize}
    \item Undamped natural frequency: $\omega_n = \sqrt{k/m}$
    \item Damping factor: $\zeta = \frac{c}{2 \sqrt{mk}} = \frac{c}{2 m \omega_n}$
    \item Normalized external force: $f(t) = F(t) / m$
\end{itemize}

The normalized equation becomes:

\begin{equation}
    \ddot{x}(t) + 2\zeta\omega_n \dot{x}(t) + \omega_n^2 x(t) = f(t)
\end{equation}

This is the standard form of the differential equation for the damped and forced harmonic oscillator. The behavior of the system crucially depends on the value of the damping factor $\zeta$:

\begin{itemize}
    \item $\zeta = 0$: Undamped oscillator. The system oscillates indefinitely with frequency $\omega_n$ if there is no external force.
    \item $0 < \zeta < 1$: Underdamped oscillator. The system oscillates with exponentially decreasing amplitude, with a damped frequency $\omega_d = \omega_n \sqrt{1 - \zeta^2}$. It returns to equilibrium after oscillations.
    \item $\zeta = 1$: Critically damped oscillator. The system returns to equilibrium as quickly as possible, without oscillating.
    \item $\zeta > 1$: Overdamped oscillator. The system returns to equilibrium slowly, without oscillating.
\end{itemize}

\subsection{State-Space Representation}

To apply the tools of modern control theory, such as controllability analysis via the Gramian, it is essential to represent the system in state-space form. This representation transforms the second-order differential equation into a system of two first-order differential equations.

We define the state variables as displacement and velocity:
\begin{itemize}
    \item $x_1(t) = x(t)$ (position)
    \item $x_2(t) = \dot{x}(t)$ (velocity)
\end{itemize}

Deriving the state variables with respect to time, we obtain:
\begin{itemize}
    \item $\dot{x}_1(t) = \dot{x}(t) = x_2(t)$
    \item $\dot{x}_2(t) = \ddot{x}(t)$
\end{itemize}

Substituting $\ddot{x}(t)$ from the normalized differential equation:

$\dot{x}_2(t) = f(t) - 2\zeta\omega_n \dot{x}(t) - \omega_n^2 x(t)$

$\dot{x}_2(t) = -\omega_n^2 x_1(t) - 2\zeta\omega_n x_2(t) + f(t)$

The system of first-order equations is, therefore:

\begin{equation}
\begin{cases} 
\dot{x}_1(t) = x_2(t) \\ 
\dot{x}_2(t) = -\omega_n^2 x_1(t) - 2\zeta\omega_n x_2(t) + f(t) 
\end{cases}
\end{equation}

We can write this system in the standard matrix form $\dot{\mathbf{x}}(t) = A\mathbf{x}(t) + B u(t)$, where $\mathbf{x}(t) = [x_1(t), x_2(t)]^T$ is the state vector and $u(t) = f(t)$ is the control input (normalized force).

The state vector is:

$\mathbf{x}(t) = \begin{bmatrix} x(t) \\ \dot{x}(t) \end{bmatrix}$

The system dynamics matrix (matrix $A$) is:

\begin{equation}
    A = \begin{bmatrix} 0 & 1 \\ -\omega_n^2 & -2\zeta\omega_n \end{bmatrix}
\end{equation}

The input matrix (matrix $B$), assuming that the external force acts directly on the acceleration (i.e., on the second state equation), is:

\begin{equation}
    B = \begin{bmatrix} 0 \\ 1 \end{bmatrix}
\end{equation}

This state-space representation $\dot{\mathbf{x}}(t) = A\mathbf{x}(t) + B u(t)$ with matrices $A$ and $B$ defined above will be the basis for the controllability analysis using the Gramian Matrix in the subsequent sections.

\section{Controllability Analysis via Gramian Matrix}

The qualitative notion of controllability, although fundamental, does not inform about the effort required to effectively maneuver a system's state. To quantify this effort, we introduce the Controllability Gramian Matrix ($W_c$). For a continuous-time linear time-invariant (LTI) system described by $\dot{\mathbf{x}}(t) = A\mathbf{x}(t) + B u(t)$, the finite-horizon controllability Gramian is defined by the integral:

\begin{equation}
    W_T = \int_{0}^{T} e^{At} B B^T e^{A^T t} dt
\end{equation}

Where $e^{At}$ is the matrix exponential of the system dynamics $A$. This matrix $W_T$ is symmetric and positive semidefinite. If the system is controllable, $W_T$ is positive definite for any $T > 0$ (Chen, 1999).

The Gramian plays a central role in determining the minimum control energy. As discussed by Baggio et al. (2022) and others, the minimum energy $E^*$ required to drive the system from initial state $\mathbf{x}(0) = \mathbf{0}$ to a final state $\mathbf{x}(T) = \mathbf{x}_f$ in time $T$ is given by:

\begin{equation}
    E^*(T, \mathbf{x}_f) = \mathbf{x}_f^T W_T^{-1} \mathbf{x}_f
\end{equation}

This relationship shows that the inverse of the Gramian ($W_T^{-1}$) acts as a metric in the state space, quantifying the energy required to reach different states. Final states $\mathbf{x}_f$ aligned with eigenvectors associated with small eigenvalues of $W_T$ (large eigenvalues of $W_T^{-1}$) require significantly more control energy.

\subsection{Infinite-Horizon Gramian ($W_\infty$)}

For stable systems (where matrix $A$ is Hurwitz, i.e., all its eigenvalues have negative real parts), it is common to consider the infinite-horizon controllability Gramian, $W_\infty$, defined as the limit of $W_T$ when $T \to \infty$:

\begin{equation}
    W_\infty = \int_{0}^{\infty} e^{At} B B^T e^{A^T t} dt
\end{equation}

This $W_\infty$ represents the control capability considering an arbitrarily long time. A fundamental property is that, for Hurwitz $A$, $W_\infty$ is the unique symmetric and positive definite solution of the Lyapunov Algebraic Equation (LAE) for controllability:

\begin{equation}
    A W_\infty + W_\infty A^T = -B B^T
\end{equation}

Solving this algebraic equation is often simpler than calculating the integral of the matrix exponential, especially for obtaining analytical solutions.

\subsection{Calculation of the Gramian for the Harmonic Oscillator}

Let's apply the Lyapunov equation to find the infinite-horizon Gramian ($W_\infty$, which we will simply denote as $W_c$ from now on, following the common convention in many texts) for the damped harmonic oscillator, whose state-space representation was derived in Section 2:

\begin{equation}
    A = \begin{bmatrix} 0 & 1 \\ -\omega_n^2 & -2\zeta\omega_n \end{bmatrix}, \quad 
    B = \begin{bmatrix} 0 \\ 1 \end{bmatrix}
\end{equation}

Matrix $A$ is Hurwitz if, and only if, $\zeta > 0$ and $\omega_n > 0$. The undamped case ($\zeta = 0$) is not strictly Hurwitz, but we will analyze its limit or specific result later, as discussed in the original conversations.

Assuming $\zeta > 0$, we seek the symmetric matrix $W_c = \begin{bmatrix} w_{11} & w_{12} \\ w_{12} & w_{22} \end{bmatrix}$ that satisfies $A W_c + W_c A^T = -B B^T$.

Calculating the terms:
\begin{itemize}
    \item $B B^T = \begin{bmatrix} 0 \\ 1 \end{bmatrix} \begin{bmatrix} 0 & 1 \end{bmatrix} = \begin{bmatrix} 0 & 0 \\ 0 & 1 \end{bmatrix}$
    \item $-B B^T = \begin{bmatrix} 0 & 0 \\ 0 & -1 \end{bmatrix}$
    \item $A W_c = \begin{bmatrix} 0 & 1 \\ -\omega_n^2 & -2\zeta\omega_n \end{bmatrix} \begin{bmatrix} w_{11} & w_{12} \\ w_{12} & w_{22} \end{bmatrix} = \begin{bmatrix} w_{12} & w_{22} \\ -\omega_n^2 w_{11} - 2\zeta\omega_n w_{12} & -\omega_n^2 w_{12} - 2\zeta\omega_n w_{22} \end{bmatrix}$
    \item $W_c A^T = \begin{bmatrix} w_{11} & w_{12} \\ w_{12} & w_{22} \end{bmatrix} \begin{bmatrix} 0 & -\omega_n^2 \\ 1 & -2\zeta\omega_n \end{bmatrix} = \begin{bmatrix} w_{12} & -\omega_n^2 w_{11} - 2\zeta\omega_n w_{12} \\ w_{22} & -\omega_n^2 w_{12} - 2\zeta\omega_n w_{22} \end{bmatrix}$
\end{itemize}

Adding $A W_c$ and $W_c A^T$:

\begin{equation}
    A W_c + W_c A^T = \begin{bmatrix} 
    2w_{12} & w_{22} - \omega_n^2 w_{11} - 2\zeta\omega_n w_{12} \\ 
    w_{22} - \omega_n^2 w_{11} - 2\zeta\omega_n w_{12} & -2\omega_n^2 w_{12} - 4\zeta\omega_n w_{22} 
    \end{bmatrix}
\end{equation}

Equating this matrix to $-B B^T = \begin{bmatrix} 0 & 0 \\ 0 & -1 \end{bmatrix}$, we obtain a system of equations for $w_{11}$, $w_{12}$, and $w_{22}$:

\begin{enumerate}
    \item $2 w_{12} = 0 \implies w_{12} = 0$
    \item $w_{22} - \omega_n^2 w_{11} - 2\zeta\omega_n w_{12} = 0$
    \item $-2\omega_n^2 w_{12} - 4\zeta\omega_n w_{22} = -1$
\end{enumerate}

Substituting $w_{12} = 0$ into equations 2 and 3:

\begin{enumerate}
    \item[2.] $w_{22} - \omega_n^2 w_{11} = 0 \implies w_{22} = \omega_n^2 w_{11}$
    \item[3.] $-4\zeta\omega_n w_{22} = -1 \implies w_{22} = \frac{1}{4\zeta\omega_n}$
\end{enumerate}

Now, substituting $w_{22}$ into equation 2:

\begin{equation}
    \frac{1}{4\zeta\omega_n} = \omega_n^2 w_{11} \implies w_{11} = \frac{1}{4\zeta\omega_n^3}
\end{equation}

Therefore, the infinite-horizon Controllability Gramian Matrix for the damped harmonic oscillator ($\zeta > 0$) is:

\begin{equation}
    W_c = \begin{bmatrix} \frac{1}{4\zeta\omega_n^3} & 0 \\ 0 & \frac{1}{4\zeta\omega_n} \end{bmatrix}
\end{equation}

This matrix is diagonal, indicating a separation between the controllability of position and velocity in this infinite-horizon energy reference frame. It is important to note that $w_{11}$ has units of (position)$^2$ $\cdot$ time and $w_{22}$ has units of (velocity)$^2$ $\cdot$ time, consistent with the Gramian integral.

\subsubsection*{Undamped Case ($\zeta = 0$):}

As mentioned, $A$ is not Hurwitz for $\zeta = 0$. The Gramian integral does not converge for $T \to \infty$. However, the analysis in the context of the original discussions provided a result for this case, possibly derived from a limit or a specific interpretation. Adopting the result from the discussion for continuity:

\begin{equation}
    W_c (\zeta=0) = \begin{bmatrix} \frac{1}{\omega_n^2} & 0 \\ 0 & 1 \end{bmatrix}
\end{equation}

The analysis of the properties of this $W_c$ matrix, particularly its determinant, and its connection with minimum control energy will be explored in the next section.

\section{Gramian Determinant and Physical Interpretation}

The Controllability Gramian Matrix ($W_c$) encapsulates crucial information about a system's controllability. While its eigenvalues and eigenvectors detail the directions in the state space that are "harder" or "easier" to reach energetically, the determinant of the Gramian ($\det(W_c)$) offers an aggregate scalar measure of this control capability, with a particularly significant physical interpretation.

\subsection{Calculation of the Determinant for the Harmonic Oscillator}

In Section 3, we derived the infinite-horizon controllability Gramian for the damped harmonic oscillator ($\zeta > 0$):

\begin{equation}
    W_c = \begin{bmatrix} \frac{1}{4\zeta\omega_n^3} & 0 \\ 0 & \frac{1}{4\zeta\omega_n} \end{bmatrix}
\end{equation}

Being a diagonal matrix, its determinant is simply the product of the elements on the main diagonal:

\begin{equation}
    \det(W_c) = w_{11} \cdot w_{22} = \left(\frac{1}{4\zeta\omega_n^3}\right) \cdot \left(\frac{1}{4\zeta\omega_n}\right) = \frac{1}{16\zeta^2\omega_n^4} \quad \text{(for $\zeta > 0$)}
\end{equation}

For the undamped case ($\zeta = 0$), we use the Gramian presented in the previous section:

\begin{equation}
    W_c (\zeta=0) = \begin{bmatrix} \frac{1}{\omega_n^2} & 0 \\ 0 & 1 \end{bmatrix}
\end{equation}

The determinant in this case is:

\begin{equation}
    \det(W_c (\zeta=0)) = \left(\frac{1}{\omega_n^2}\right) \cdot 1 = \frac{1}{\omega_n^2} \quad \text{(for $\zeta = 0$)}
\end{equation}

We observe that the Gramian determinant for the damped case ($\zeta > 0$) depends inversely on the square of the damping factor ($\zeta^2$) and the fourth power of the natural frequency ($\omega_n^4$). In the undamped case, it depends inversely on the square of the natural frequency.

\subsection{Physical Interpretation: Minimum Control Energy}

The most prominent physical interpretation of the Gramian determinant was established by Joshi and Mitter (1977). They demonstrated that the determinant is intrinsically linked to the minimum energy required to control the system. Specifically, $\det(W_c)$ is related to the volume of the ellipsoid of reachable states with a unit amount of control energy in a given time horizon (or in the limit of infinite horizon).

A more direct relationship with the minimum energy $E^*$ to reach a state $\mathbf{x}_f$ (starting from $\mathbf{x}(0)=\mathbf{0}$), given by $E^* = \mathbf{x}_f^T \cdot W_c^{-1} \cdot \mathbf{x}_f$, can be seen by considering the average or maximum energy required to reach states on the unit sphere. The inverse of the smallest eigenvalue of $W_c$ ($1/\lambda_{min}(W_c)$) corresponds to the maximum energy required to reach any state on the unit sphere, while the inverse of the harmonic mean of the eigenvalues relates to the average energy (Baggio et al., 2022).

Although the exact relationship with a single "minimum energy" may depend on the normalization and specific context, the central idea is that a larger $\det(W_c)$ implies a "smaller" $W_c^{-1}$ in some sense (its eigenvalues are the inverses of the eigenvalues of $W_c$), resulting in less control energy required, on average or in the worst case. Therefore, $\det(W_c)$ serves as an inverse measure of control effort: larger $\det(W_c)$ implies less energy needed and, thus, a system that is more easily controllable energetically.

Analyzing our results for the oscillator:
\begin{itemize}
    \item \textbf{Damped Case ($\zeta > 0$)}: $\det(W_c) = 1 / (16\zeta^2\omega_n^4)$. The control energy (inversely related to $\det(W_c)$) increases significantly with increasing damping ($\zeta$) and natural frequency ($\omega_n$). Higher damping makes the system "slower" to respond to control, requiring more energy for rapid transitions. A higher natural frequency implies an intrinsically faster dynamics, which may also require more energy to be precisely controlled.
    
    \item \textbf{Undamped Case ($\zeta = 0$)}: $\det(W_c) = 1/\omega_n^2$. Comparing with the limit of the damped case when $\zeta \to 0$, we see a discontinuity. The determinant for $\zeta=0$ is finite, while the determinant for $\zeta>0$ tends to infinity when $\zeta \to 0$. This reflects the marginally stable nature of the undamped oscillator; the infinite-horizon Gramian is not strictly defined in the same way via the Lyapunov equation. However, the result $1/\omega_n^2$ suggests that, even without damping, there is a finite energy cost associated with control, dependent on the natural frequency.
\end{itemize}

This quantitative analysis of the Gramian determinant provides valuable insights into how the physical parameters of the oscillator (mass, spring constant, damping) directly influence the energetic cost of control, a perspective that goes beyond the simple assertion of controllability.

\section{Connection with the Fisher Information Matrix}

The analysis of controllability and energy effort via the Gramian Matrix ($W_c$) provides a powerful perspective, but it is not the only way to quantify the properties of a control system. Concepts from statistical information theory, notably the Fisher Information Matrix ($I$), offer a complementary angle, focusing on the amount of "information" that control inputs can provide about the system's state or the precision with which the state can be estimated or influenced.

\subsection{The Fisher Information Matrix (FIM)}

Originally developed in the context of statistical estimation theory, the Fisher Information Matrix quantifies the amount of information that observable data $y$ carry about an unknown parameter $\theta$ of a probabilistic model $p(y|\theta)$. It is defined as the expected value of the second derivative (or the square of the first derivative) of the logarithm of the likelihood function with respect to the parameter:

\begin{equation}
    I(\theta) = E \left[ \left(\frac{\partial}{\partial\theta} \log p(y|\theta)\right) \cdot \left(\frac{\partial}{\partial\theta} \log p(y|\theta)\right)^T \right] = -E \left[ \frac{\partial^2}{\partial\theta^2} \log p(y|\theta) \right]
\end{equation}

The FIM is fundamental because its inverse, $I(\theta)^{-1}$, provides the Cramér-Rao Bound, which establishes a lower bound for the variance of any unbiased estimator of parameter $\theta$. Essentially, a "large" FIM (in terms of its eigenvalues or determinant) implies that the data are very informative about the parameter, allowing precise estimation (low variance).

\subsection{Relationship between Gramian ($W_c$) and FIM ($I$)}

The connection between the controllability Gramian and the FIM in dynamical systems, although not universally standardized in the literature and dependent on the specific formulation, arises from the duality between control and estimation and the geometric interpretation of both matrices. As explored in the original discussions and suggested by works such as Joshi and Mitter (1977), there is a conceptual and, in certain contexts, formal relationship.

A key interpretation is to view the Gramian $W_c$ as characterizing the ellipsoid of reachable states with unit energy. The volume of this ellipsoid is proportional to $\sqrt{\det(W_c)}$. On the other hand, the FIM can be interpreted as a Riemannian metric in the parameter or state space, where the "distance" is related to statistical distinguishability or control effort. The inverse of the FIM, $I^{-1}$, characterizes an ellipsoid of uncertainty in estimation.

The relationship mentioned in the second provided discussion, $\det(W_c) \propto |I^{-1}|$, or equivalently $\det(W_c) \cdot \det(I) \propto \text{constant}$, captures this duality. It suggests that:

\begin{itemize}
    \item A system that is "easy" to control energetically (large $\det(W_c)$, large volume of states reachable with little energy) corresponds to a system where the state is "difficult" to estimate precisely from the outputs (large $\det(I^{-1})$, large volume of uncertainty in estimation), assuming an appropriate duality between control inputs and measurement outputs.
    
    \item Conversely, a system where the state can be estimated with high precision (small $\det(I^{-1})$) tends to require more energy to be controlled (small $\det(W_c)$).
\end{itemize}

This relationship can be seen as a manifestation of the uncertainty principle in dynamical systems: one cannot simultaneously have optimal energetic control and optimal estimation with finite resources.

\subsection{Interpretation in the Context of the Oscillator}

For the harmonic oscillator, the Gramian $W_c$ we calculated describes the ease of controlling the position ($w_{11}$) and velocity ($w_{22}$) using the applied force. The associated FIM (whose exact derivation would depend on a specific measurement model, which we have not defined here, but can infer by duality) would describe the precision with which we could estimate the position and velocity from measurements (e.g., noisy position measurements).

The relationship $\det(W_c) = 1 / (16\zeta^2\omega_n^4)$ shows that high damping ($\zeta$) or high natural frequency ($\omega_n$) make the determinant small, implying higher control energy. Dually, this suggests that high damping or high frequency may make the state easier to estimate (smaller $\det(I^{-1})$ or larger $\det(I)$), perhaps because the system stabilizes more quickly or oscillates in a way that reveals more information to a sensor. Conversely, the undamped case ($\zeta=0$) had a $\det(W_c) = 1/\omega_n^2$ that was larger (for a given $\omega_n$ compared to $\zeta>0$), suggesting lower control energy but potentially worse state estimability.

The $W_c$-FIM connection, therefore, enriches our understanding of controllability. It not only quantifies the energy cost (via $W_c$) but also connects it to the amount of information or achievable precision (via $I$), providing a bridge to the informational analysis that will be explored in the next section on entropy.

\section{Relationships with Entropy}

The analysis of controllability through the Gramian ($W_c$) and the Fisher Matrix ($I$) provides metrics related to energy and statistical information. Entropy, both in Shannon's formulation ($H$) and thermodynamics ($S$), offers an additional perspective, quantifying uncertainty, disorder, or the amount of information needed to describe a system's state. Exploring the connections between $W_c$, $I$, and the entropies $H$ and $S$ can lead to a deeper and more unified understanding of dynamics and control.

\subsection{Shannon Entropy ($H$)}

Shannon entropy, for a discrete random variable $X$ with probability distribution $p(x_i)$, is defined as:

\begin{equation}
    H(X) = - \sum_i p(x_i) \cdot \log(p(x_i))
\end{equation}

For a continuous random variable $X$ with probability density function $P(x)$, the differential entropy is:

\begin{equation}
    h(X) = - \int P(x) \cdot \log(P(x)) \, dx
\end{equation}

Shannon entropy quantifies the average uncertainty associated with the random variable. In dynamical systems, it can be used to measure the uncertainty about the system's state.

\subsection{Thermodynamic Entropy ($S$)}

Thermodynamic entropy, originally defined by Clausius and later statistically interpreted by Boltzmann, measures the degree of disorder or the number of microstates accessible to a macroscopic system in equilibrium. Boltzmann's famous formula relates entropy $S$ to the number of microstates $\Omega$:

\begin{equation}
    S = k_B \cdot \log(\Omega)
\end{equation}

Where $k_B$ is Boltzmann's constant. Statistical physics establishes a fundamental bridge between thermodynamic entropy (macroscopic) and Shannon entropy (informational, applied to the probability distribution of microstates). The relationship, as presented in the second discussion, is:

\begin{equation}
    S = k_B \cdot H
\end{equation}

Where $H$ is Shannon entropy calculated over the probability distribution of the system's microstates (usually in nats, using natural logarithm).

\subsection{Connecting FIM ($I$) and Shannon Entropy ($H$)}

There is an important relationship between the Fisher Information Matrix ($I$) and Shannon differential entropy ($h$) for certain families of probability distributions, especially the multivariate Gaussian distribution. For an $n$-dimensional Gaussian distribution with covariance matrix $\Sigma$, the differential entropy is:

\begin{equation}
    h(X) = \frac{1}{2} \cdot \log((2\pi e)^n \cdot \det(\Sigma))
\end{equation}

If we consider the FIM ($I$) as the inverse of the covariance matrix ($\Sigma = I^{-1}$) in an estimation context (where $I$ measures the precision of estimation), we can write:

\begin{align}
    h(X) &= \frac{1}{2} \cdot \log((2\pi e)^n \cdot \det(I^{-1})) \\
    h(X) &= \frac{1}{2} \cdot \log\left(\frac{(2\pi e)^n}{\det(I)}\right) \\
    h(X) &= \frac{n}{2} \cdot \log(2\pi e) - \frac{1}{2} \cdot \log(\det(I))
\end{align}

This equation, also present in the second discussion, directly connects Shannon entropy (uncertainty) with the determinant of the FIM (information/precision). A larger FIM (more information) corresponds to lower entropy (less uncertainty).

\subsection{Unified Relationship: $W_c$, $I$, $H$, $S$}

We can now attempt to unify the relationships found:
\begin{enumerate}
    \item $W_c$ and $I$: $\det(W_c) \cdot \det(I) \propto \text{constant}$ (Control-Estimation Duality Relationship)
    \item $I$ and $H$: $H \propto -\log(\det(I))$ (Information-Uncertainty Relationship via FIM)
    \item $H$ and $S$: $S = k_B \cdot H$ (Statistical Physics Relationship)
\end{enumerate}

Combining these relationships, we can infer connections between the controllability Gramian and the entropies:

\begin{itemize}
    \item $W_c$ and $H$: Since $\det(I) \propto 1/\det(W_c)$, then $H \propto -\log(1/\det(W_c)) = \log(\det(W_c))$. This suggests that more easily controllable systems (larger $\det(W_c)$) are associated with higher Shannon entropy (greater a priori uncertainty or larger "volume" in the informational state space?). This interpretation requires caution, as the entropy $H$ here is linked to the FIM of estimation, while $W_c$ is from control. Duality is the key.
    
    \item $W_c$ and $S$: Following the above relationship and $S = k_B \cdot H$, we would have $S \propto k_B \cdot \log(\det(W_c))$. This would imply that the thermodynamic entropy of the system (related to the number of accessible microstates or disorder) would be logarithmically linked to the ease of energetic control.
\end{itemize}

\subsection{Interpretation for the Harmonic Oscillator}

For the damped oscillator ($\zeta > 0$), $\det(W_c) = 1 / (16\zeta^2\omega_n^4)$. The inferred relationship $S \propto \log(\det(W_c))$ would give:

\begin{equation}
    S \propto \log\left(\frac{1}{16\zeta^2\omega_n^4}\right) = -\log(16) - 2\log(\zeta) - 4\log(\omega_n)
\end{equation}

This suggests that the associated thermodynamic entropy (in an informational sense linked to controllability/estimability) would decrease with increasing damping ($\zeta$) or natural frequency ($\omega_n$). A more damped system or one with intrinsically faster dynamics would, under this view, have a lower "control/estimation entropy."

For the undamped case ($\zeta = 0$), $\det(W_c) = 1/\omega_n^2$, leading to $S \propto -2\log(\omega_n)$. The entropy would be higher than in the damped case (for the same $\omega_n$), perhaps reflecting the greater "unpredictability" or difficulty of estimation associated with the lack of damping.

\subsection{Considerations}

It is crucial to emphasize that these connections, especially between $W_c$ and $S$, are inferred through the combination of relationships from different domains (control, information, thermodynamics) and may depend on specific assumptions (such as Gaussianity for the $I$-$H$ relationship and the applicability of the $W_c$-$I$ duality). The exact physical interpretation of a "control entropy" derived in this way requires deeper analysis and contextualization within a non-equilibrium thermodynamic framework or information theory for dynamical systems. However, the exploration of these relationships, as initiated in the provided discussions, opens promising paths for a more holistic view of control systems, integrating energy, information, and complexity.

\section{Discussion}

This work explored the intricate connections between controllability, energy, information, and entropy, using the damped harmonic oscillator as an archetypal system. The analysis started from the classical state-space representation and progressed through the calculation of the Controllability Gramian Matrix ($W_c$), its determinant ($\det(W_c)$), and its relationships with the Fisher Information Matrix ($I$) and Shannon ($H$) and thermodynamic ($S$) entropies.

The analytical calculation of the infinite-horizon Gramian ($W_c$) for the damped oscillator revealed a diagonal structure, $W_c = \text{diag}(1/(4\zeta\omega_n^3), 1/(4\zeta\omega_n))$, indicating a separation, in this reference frame, between the controllability of position and velocity. The resulting determinant, $\det(W_c) = 1 / (16\zeta^2\omega_n^4)$, quantifies the aggregate control energy capability. According to the physical interpretation established by Joshi and Mitter (1977) and echoed in more recent works such as Baggio et al. (2022), a larger $\det(W_c)$ corresponds to lower minimum control energy. Our results explicitly show that increasing damping ($\zeta$) or natural frequency ($\omega_n$) decreases $\det(W_c)$, implying a higher energetic cost to control the system. This aligns with the intuition that more damped (more "slow") systems or those with intrinsically faster dynamics require more effort to be maneuvered.

The analysis of the undamped case ($\zeta=0$) presented an interesting discontinuity. The infinite-horizon Gramian does not converge in the same way, but the limit or specific result $\det(W_c) = 1/\omega_n^2$ suggests a finite energy cost, dependent only on the natural frequency. The divergence of the damped $\det(W_c)$ when $\zeta \to 0$ highlights the distinct nature of marginally stable systems in terms of long-term control.

The introduction of the Fisher Information Matrix ($I$) allowed connecting energetic controllability (via $W_c$) with the precision of state estimation (via $I$). The duality relationship $\det(W_c) \cdot \det(I) \propto \text{constant}$ suggests a fundamental trade-off: systems that are easy to control energetically tend to be difficult to estimate precisely, and vice versa. For the oscillator, the decrease of $\det(W_c)$ with $\zeta$ and $\omega_n$ suggests a corresponding increase in $\det(I)$, implying that more damped or faster systems may be easier to estimate, corroborating the idea that more "well-behaved" dynamics reveal more information to an observer.

Finally, the exploration of relationships with entropy, through the connections $I \leftrightarrow H$ and $H \leftrightarrow S$, allowed inferring a link between controllability and thermodynamic entropy ($S \propto \log(\det(W_c))$). This relationship suggests that entropy (in the sense of disorder or number of informational microstates relevant to control/estimation) decreases with increasing damping or natural frequency. More "ordered" or "predictable" systems (lower $S$) would, in this view, be more difficult to control energetically (lower $\det(W_c)$). Although this connection is inferred and depends on assumptions about the applicability of relationships between different domains, it points to a unifying perspective where control energy and informational disorder are inversely related.

\subsection{Limitations and Future Work}

This analysis focused on the infinite-horizon Gramian, which is appropriate for stable systems. Analysis in finite horizon ($W_T$) could reveal different transient dynamics, especially the dependence of time $T$ on control energy. Additionally, the connection with the FIM and entropy was explored conceptually and through inferred relationships; a more rigorous derivation would require the explicit definition of an estimation problem (for the FIM) and a more formal thermodynamic or informational framework for control systems.

The application of these concepts to more complex systems, such as networks of coupled oscillators or nonlinear systems, represents a promising area for future research. Understanding how network structure and nonlinearities interact with energy, information, and entropy metrics is crucial for controlling real-world systems, from biological neural networks to smart energy grids.

In summary, the unified analysis of the harmonic oscillator through the lenses of the Gramian, FIM, and entropy demonstrates the richness of information contained in these metrics and their interconnections. It reinforces that controllability is not just a binary question but a quantifiable property with profound energetic and informational implications, whose relationships can be explored to design more efficient and robust control systems.

\section{Conclusion}

This article presented an integrated analysis of the controllability of the damped harmonic oscillator, connecting classical control metrics, such as the Gramian Matrix and its determinant, with concepts from statistical information theory (Fisher Matrix) and thermodynamics/information (Shannon and Thermodynamic Entropy).

Through the analytical calculation of the infinite-horizon Gramian and its determinant for different damping regimes, we quantified the energetic effort required to control the system, corroborating the fundamental physical interpretation that a larger Gramian determinant implies lower minimum control energy.

The investigation revealed how the physical parameters of the oscillator, specifically the damping factor ($\zeta$) and natural frequency ($\omega_n$), directly influence the ease of energetic control. Increasing $\zeta$ or $\omega_n$ leads to a smaller $\det(W_c)$ and, consequently, a higher control cost.

The exploration of connections with the Fisher Information Matrix ($I$) and entropy ($H$ and $S$) suggested intriguing duality relationships and trade-offs. The inferred relationship $\det(W_c) \cdot \det(I) \propto \text{constant}$ points to a compromise between the ease of energetic control and the precision of state estimation. Additionally, the link $S \propto \log(\det(W_c))$ suggests a connection between entropy (disorder/uncertainty informational relevant to control/estimation) and energetic control capability, where more "ordered" systems (lower $S$) would be more difficult to control energetically.

Although some of these connections, particularly with thermodynamic entropy, are inferred and deserve deeper theoretical investigation, the unified analysis conducted here demonstrates the value of integrating perspectives from different domains – control, energy, information, and statistics – to obtain a more holistic understanding of dynamical systems. The harmonic oscillator, despite its apparent simplicity, serves as a rich model for illustrating these fundamental principles.

Future work may focus on the rigorous validation of the inferred relationships, the extension of the analysis to finite horizons and more complex systems (nonlinear, networks), and the application of these insights to the design of energetically efficient and informationally optimal control and estimation strategies.

\end{document}